\documentclass[leqno,final]{amsart}
\usepackage{amssymb,amsmath,latexsym,amsfonts,graphics,color,epsf,epsfig,psfrag}

\newcommand{\st}{{\,|\,}}		%general
\renewcommand{\frak}{\mathfrak}		%general
\renewcommand{\Bbb}{\mathbb}		%general
\newcommand{\ds}{\displaystyle}		%general
\newcommand{\gdn}{\mathbb{G}_{d,n}}		%theorem.tex
\newcommand{\idn}{I(d,n)}		%theorem.tex
\newcommand{\mon}{\frak S}		%proof.tex
\renewcommand{\part}{\frak B}		%proof.tex

\newcommand{\product}{\prod}		%proof.tex
\newcommand{\path}{\Lambda}		%interpretation.tex
\newcommand{\bstart}{\beta({\mbox{\begin{rm}start\end{rm}}})}
\newcommand{\bstartj}{\beta_j({\mbox{\begin{rm}start\end{rm}}})}
\newcommand{\bfinish}{\beta({\mbox{\begin{rm}finish\end{rm}}})}
\newcommand{\bfinishj}{\beta_j({\mbox{\begin{rm}finish\end{rm}}})}
\newtheorem{theorem}{Theorem}%[section]
\newcommand{\bthm}{\begin{theorem}}
\newcommand{\ethm}{\end{theorem}}

\newtheorem{proposition}[theorem]{Proposition}
\newcommand{\bpr}{\begin{proposition}}
\newcommand{\bprop}{\begin{proposition}}
\newcommand{\epr}{\end{proposition}}
\newcommand{\eprop}{\end{proposition}}

\newtheorem{definition}[theorem]{Definition}
\newcommand{\bdefn}{\begin{definition}\begin{rm}}
\newcommand{\edefn}{\end{rm}\end{definition}}

\newtheorem{example}{Example}
\newcommand{\bexample}{\begin{example}\begin{rm}}
\newcommand{\eexample}{\ $\Box$\end{rm}\end{example}}

\newcommand{\bexampletwo}{\begin{example}\begin{rm}}
\newcommand{\eexampletwo}{\end{rm}\end{example}}

\newtheorem{remark}[theorem]{Remark}
\newcommand{\bremark}{\begin{remark}\begin{rm}}
\newcommand{\eremark}{\end{rm}\end{remark}}

\newtheorem{corollary}[theorem]{Corollary}
\newcommand{\bcor}{\begin{corollary}}
\newcommand{\ecor}{\end{corollary}}

\newtheorem{exercise}[theorem]{Exercise}
\newcommand{\bex}{\begin{exercise}\begin{rm}}
\newcommand{\eex}{\end{rm}\end{exercise}}

\newtheorem{lemma}[theorem]{Lemma}
\newcommand{\blem}{\begin{lemma}}
\newcommand{\blemma}{\begin{lemma}}
\newcommand{\elemma}{\end{lemma}}
\newcommand{\elem}{\end{lemma}}

\newenvironment{proofone}[1][]{{\sc Proof #1}:\ }{\hfill$\Box$}
\newcommand{\bpf}{\begin{proofone}}
\newcommand{\bproof}{\begin{proofone}}
\newcommand{\epf}{\end{proofone}}
\newcommand{\eproof}{\end{proofone}}

\newenvironment{prooftwo}[1][]{{\sc Proof #1}:\ }{}
\newcommand{\bpftwo}{\begin{prooftwo}}
\newcommand{\bprooftwo}{\begin{prooftwo}}
\newcommand{\epftwo}{\end{prooftwo}}
\newcommand{\eprooftwo}{\end{prooftwo}}

\newenvironment{solution}{{\sc Solution}:\ }{\ $\Box$}
\newcommand{\bsol}{\begin{solution}}
\newcommand{\esol}{\end{solution}}
\newcommand{\coh}{\mathcal H}
\newcommand{\base}{\mathcal S}
\newcommand{\bring}{\base}
\newcommand{\point}{{\rm Spec}(\mathbb{C})}
\newcommand{\basis}{{\cal B}}
\newcommand{\res}{{\rm Res}}
\newcommand{\VDM}{\mathfrak V}          %giambelli
\newcommand{\Num}{\mathfrak{Num}}               %giambelli
\newcommand{\coloneq}{:=}               %general
\newcommand{\roots}{\frak R }           %restrproof.tex
\newcommand{\pos}{\frak N }           %restrproof.tex
\newcommand{\srsc}{{\cal C}_u^v}	   	%restrproof.tex
\newcommand{\srfr}{R}			%restrproof.tex
\newcommand{\face}{f}			%restrproof.tex
\newcommand{\seq}{\underline{\ell}}	%restrproof-2.tex
\newcommand{\srscprime}{{\cal C}_{u'}^{v'}}	   	%restrproof-2.tex
\newcommand{\codim}{\mbox{{\rm codim}$\,$}}

\newcommand{\lr}{\mbox{$\longrightarrow$}}
\newcommand{\bE}{E}%{\mathbb{E}}
\newcommand{\bb}{B}%{\mathbb{B}}

\newcommand{\bc}{\mathbb{C}}
\newcommand{\bp}{\mathbb{P}}

\renewcommand{\phi}{\varphi}

\newcommand{\cal}{\mathcal}

\newcommand{\ignore}[1]{}
\begin{document}
\title[Equivariant Giambelli and restriction formulas]{Equivariant Giambelli and determinantal restriction formulas
for the Grassmannian}
\author[V. Lakshmibai]{V. Lakshmibai${}^{\dag}$}
\address{Department of Mathematics\\ Northeastern University\\ Boston, MA
02115} \email{lakshmibai@neu.edu}
\thanks{${}^{\dag}$ Partially supported by NSF grant
DMS~0400679
and NSA-MDA~904-03-1-0034},
\author{K.~N.~Raghavan}
\address{Institute of Mathematical Sciences,
C.I.T.~Campus, 
Chennai 600\thinspace113 INDIA}
\email{knr@imsc.res.in}
\author{P.~Sankaran}
\address{Institute of Mathematical Sciences,
C.I.T.~Campus, 
Chennai 600\thinspace113 INDIA}
\email{sankaran@imsc.res.in}
\date{25~Aug~2005\\
\indent 2000 {\em Mathematics 
Subject Classification:} 14M15}
\begin{abstract}The main result of the paper is a determinantal formula for the
restriction to a torus fixed point of the equivariant class
of a Schubert subvariety in the torus equivariant integral cohomology
ring of the Grassmannian.    As a corollary,  we obtain an
equivariant version of the Giambelli formula.
\end{abstract}

\maketitle
The (torus) equivariant cohomology rings 
of flag varieties in general and
of the Grassmannian in particular have recently attracted
much interest.   
Here we consider the equivariant integral cohomology ring %---let us
of the Grassmannian.
Just as the ordinary Schubert classes 
form a module basis over the ordinary cohomology ring of a point
(namely the ring of integers) for the ordinary integral cohomology ring of %$X$,
the Grassmannian, 
so do the equivariant Schubert classes form a basis 
over the equivariant
cohomology of a point (namely the ordinary cohomology ring
of the classifying space of the torus) for the equivariant cohomology ring (this is true for any generalized flag variety of any type, not just the
Grassmannian).    
Again  as in
the ordinary case,  computing the structure constants of the multiplication 
with respect to this basis is an interesting problem that goes by the name
of Schubert calculus.
There is a forgetful functor from 
equivariant cohomology to ordinary cohomology so that results about the former
specialize to those about the latter.

Knutson-Tao-Woodward~\cite{ktw} and Knutson-Tao~\cite{kt} show that the structure constants, both
ordinary and equivariant, count 
solutions to certain jigsaw puzzles, thereby showing that they are
``manifestly'' positive.
In the present paper we take a very different route to computing the
equivariant structure constants.    Namely,  we try to extend
to the equivariant
case the classical approach by means of the Pieri
and Giambelli formulas.   Recall, from \cite[Eq.(10), p.146]{fulton:yt} for example,
that the Giambelli formula expresses an
arbitrary Schubert class as a polynomial with integral coefficients in certain
``special'' Schubert classes---the Chern classes of the tautological
quotient bundle---and that the Pieri formula expresses as a
linear combination of the Schubert classes the product of a special Schubert
class with an arbitrary Schubert class.    Together they can be used to
compute the structure constants.     

We only partially succeed in our attempt:  the first of the three theorems
of this paper---see \S\ref{sgiambelli} below---is an equivariant Giambelli
formula that specializes to the ordinary Giambelli formula as in 
\cite[Eq.(10), p.146]{fulton:yt}, but we still
do not have a satisfactory equivariant Pieri formula---see, however,
\S\ref{spieri} below.   The derivation in
Fulton \cite[\S14.3]{fulton:it} of the Giambelli formula
can perhaps be extended to the equivariant case,
but this is not what we do.     Instead, we deduce
the Giambelli formula 
from our second theorem which gives a certain
closed-form determinantal formula for the restriction to a torus fixed
point of an equivariant Schubert class.  

This ``restriction formula'' (Theorem~\ref{trestrdet} in \S\ref{srestrdet} below) is the point of this
paper---more so than the Giambelli, for among other things it might also
hold the key to the Pieri. 
It in turn is deduced
from Theorem~\ref{tgrobrestr}
which can be paraphrased thus:
a Gr\"obner degeneration
of an open piece around a torus fixed point
of a torus stable subvariety
computes the restriction to the fixed point
of the equivariant cohomology class of the subvariety.
Recall that a Gr\"obner degeneration comes from a Gr\"obner basis.  
It is a $1$-parameter flat degeneration. See~\S\ref{sgrobrestr} for
the precise meaning.
Such Gr\"obner degenerations at torus fixed points of
Schubert subvarieties in the Grassmannian have recently been
obtained~\cite{kl,kr,k,kl:two}.   Combining this result about degenerations
with Theorem~\ref{tgrobrestr} yields
a proof of Theorem~\ref{trestrdet}.

As pointed out to us by the referee,
Theorem~\ref{tgrobrestr} is well known.
The precise references (as indicated by the referee) are given
in~\S\ref{sgrobrestr}.  
%Nevertheless we think it makes sense
%to include a proof here for the sake of completeness and convenience.
%
%The proof in~\S\ref{sgrobrestr} of Theorem~\ref{tgrobrestr}
%uses only standard facts about
%equivariant cohomology (for which we refer to~\cite{kt} and the references therein) and about flat degenerations (for which we refer to the books~\cite{ebud,fulton:it}) .
The passage from Theorem~\ref{tgrobrestr} plus the result about degenerations to
the restriction formula involves only an elementary combinatorial 
inductive argument.  (Is
there an elegant Lindstrom-Gessel-Viennot type argument for this passage?
We do not know.)
The passage from the restriction formula to the Giambelli again involves
only elementary matrix manipulations.    

{\sc Acknowledgments:}  Parts of this work were done during visits of the
first named author to {\ttfamily Chennai Mathematical Institute\/} and 
of the other two authors to {\ttfamily The Abdus Salam International
Centre for Theoretical Physics\/}.
The hospitality of the two institutions is gratefully acknowledged.
It is with great pleasure that we thank the referee for a quick and
thorough reading
of the manuscript and for the insightful comments;  thanks in particular for
indicating how Theorem~\ref{tgrobrestr} follows from results available in
the literature.

\section{The set up}\label{ssetup}
Fix once and for all two positive integers $d$ and $n$ with $d\leq n$.
Let $V$ be an $n$-dimensional complex vector space,  and $\gdn$
the Grassmannian of $d$-dimensional linear subspaces of $V$.
The defining action of the general linear group $GL(V)$ on $V$
induces an action on $\gdn$.    We are interested in the 
$T$-equivariant integral cohomology ring $\coh$ of $\gdn$,  where $T$
is a fixed maximal torus of $GL(V)$.    

We refer to \cite[\S2]{kt} and the
references in that paper for the details that we leave out in this section.

The natural map from $\gdn$ to $\point$ induces an
$\base$-algebra structure on $\coh$,  where $\base:=H^*_T(\point)$
is the $T$-equivariant integral cohomology ring of $\point$
(namely the ordinary integral cohomology ring of the classifying space of
$T$).    The $\base$-algebra $\coh$ is independent of the choice of $T$
because any two maximal tori in $GL(V)$ are conjugate.

The choice of a maximal torus $T$ amounts to the choice of an
unordered vector space basis $\basis$ of $V$:  the elements of $T$
are precisely those invertible linear transformations for which each
element of $\basis$ is an eigenvector.    Each element $b$ of $\basis$
thus defines a character $\epsilon_b$ of $T$ that sends elements of $T$
to their respective eigenvalues with respect to $b$.   
The collection $\{\epsilon_b\st b\in\basis\}$ forms
an integral basis for the group $X(T)$ of characters of $T$.  The ring $\base$
is graded isomorphic to the symmetric algebra
of the abelian group $X(T)$ with $X(T)$ living in degree $2$.
We may therefore identify $\base$ with
the polynomial ring $\mathbb{Z}[\epsilon_b\st b\in\basis]$,  where
the $\epsilon_b$ are variables in degree $2$.

Since $\gdn$ is a smooth variety on which $T$ acts algebraically with finitely 
many fixed points,   it follows 
that $\coh$ is
a free $\base$-module with basis the (equivariant) classes of the
Schubert subvarieties.     These subvarieties are defined with respect to
a fixed Borel subgroup $B$ containing $T$:   they are the closures of the
$B$-orbits in $\gdn$.     The formulas of this paper are independent of 
the choice of $B$ because
any two such Borel subgroups are conjugate by an element in the normalizer of $T$
in $GL(V)$.    

The choice of a Borel subgroup $B$ containing $T$ amounts to the choice of an
ordering on the elements of the basis $\basis$.    Let $b_1,\ldots, b_n$ be the
elements of $\basis$ thus ordered.    We have $\base=\mathbb{Z}[\epsilon_1,
\ldots, \epsilon_n]$,  where $\epsilon_j:=\epsilon_{b_j}$.  

There is a one-to-one correspondence between the $B$-orbits and the $T$-fixed
points in $\gdn$:   each $B$-orbit contains one and only one $T$-fixed point.
The $T$-fixed points are indexed by the subsets of cardinality $d$ of $\basis$.

Denote by $\idn$ the set of subsets of cardinality $d$ of $\{1,\ldots, n\}$.
We use $u,v,w,\ldots$ to denote elements of $\idn$.   For $u$ in $\idn$,
we write $u=(u_1,\ldots,u_d)$ where $u_1,\ldots,u_d$ are the elements
of $u$ arranged in increasing order: $1\leq u_1<\ldots<u_d\leq n$.

Given $u=(u_1,\ldots,u_d)$ in $\idn$,  denote by $e^u$ the 
$T$-fixed point of $\gdn$ that is the span of $\{b_{u_1},\ldots,b_{u_d}\}$, by $X(u)$
the closure of the $B$-orbit containing $e^u$,  by $[X(u)]$ the $T$-equivariant
class in $\coh$ of the Schubert subvariety $X(u)$,  and by $[X(u)]_{\rm cl}$ the
ordinary cohomology class of $X(u)$. 

For each $T$-fixed point $e^v$, $v$ in $\idn$, we have a natural ``restriction''
map $\res_v\!:\ \coh:=H^*_T(\gdn)\to\base:=H^*_T(e^v)$ induced by the inclusion
of $\{e^v\}$ in $\gdn$.    The direct product
of these is an injection of rings:
\begin{equation}\label{einjection}
\product \res_v:\ \ \coh\hookrightarrow \product\limits_{v\in\idn}H^*_T(e^v) 
\end{equation}
For $u$ and $v$ in $\idn$, denote by $[X(u)]|_v$ the image
in $\base$ under $\res_v$   
of the equivariant class $[X(u)]$.
The image of $\coh$ under $\prod\res_v$ has a neat
description
but we will have no use for this here:
%\footnote{
a tuple $(\alpha_v)_{v\in\idn}$ in $\product_{\idn}\base$ belongs to the image of 
$\coh$ under $\prod\res_v$
if and only if, whenever $w$ and $x$ in $\idn$ are so related that there exist
integers $i$ and $j$, $1\leq i, j\leq n$, with $x=(w\cup\{j\})\setminus\{i\}$,
it holds that $\epsilon_j-\epsilon_i$ divides $\alpha_x-\alpha_w$.
%}

\raggedbottom
\section{An equivariant Giambelli formula}\label{sgiambelli}
Given $u=(u_1,\ldots,u_d)$ in $I(d,n)$,  set
\[
	\lambda_1:= n-d+1-u_1,\quad\ldots\quad\lambda_i:=n-d+i-u_i,\quad\ldots\quad \lambda_d:=n-u_d. \]
Then $n-d\geq\lambda_1\geq\lambda_2\geq\ldots\geq\lambda_d\geq0$.    

If $u$ is such that $\lambda_2=\ldots=\lambda_d=0$,  we call the Schubert variety $X(u)$
and its cohomology class {\em special};   furthermore the equivariant
and ordinary cohomology classes $[X(u)]$ and $[X(u)]_{\rm cl}$ are denoted
instead by $[\lambda_1]$ and $[\lambda_1]_{\rm cl}$.
We extend the terminology and notation
to all integers  
by setting $[p]:=0$
if $p$ is outside the range $0$, $1$, \ldots, $n-d$.

Observe that $[p]$ belongs to $H_T^{2p}(\gdn)$,
which explains the notation.

The classical Giambelli formula gives an expression for an arbitrary Schubert
class in the ordinary cohomology ring of the Grassmannian $\gdn$
as the determinant of a $d\times d$ matrix whose entries are special classes.    For
$u=(u_1,\ldots,u_d)$ in $I(d,n)$,   we have,  from \cite[Eq.(10),~page146]{fulton:yt} for
example,
$[X(u)]_{\text{\rm cl}}=$ 
\[
\left|
\begin{array}{cccccc}
[\lambda_1]_{\text{\rm cl}} &
[\lambda_1+1]_{\text{\rm cl}} &
\ldots &
[\lambda_1+j-1]_{\text{\rm cl}} & \ldots &
[\lambda_1+d-1]_{\text{\rm cl}} \\[2mm]
[\lambda_2-1]_{\text{\rm cl}} &
[\lambda_2]_{\text{\rm cl}} &
\ldots &
[\lambda_2+j-2]_{\text{\rm cl}} & \ldots &
[\lambda_2+d-2]_{\text{\rm cl}} \\[2mm]
\vdots & \vdots & \vdots & \vdots & \vdots & \vdots\\[2mm]
[\lambda_i+1-i]_{\text{\rm cl}} &
[\lambda_i+2-i]_{\text{\rm cl}} &
\ldots &
[\lambda_i+j-i]_{\text{\rm cl}} & \ldots &
[\lambda_i+d-i]_{\text{\rm cl}} \\[2mm]
\vdots & \vdots & \vdots & \vdots & \vdots & \vdots\\[2mm]
[\lambda_d+1-d]_{\text{\rm cl}} &
[\lambda_d+2-d]_{\text{\rm cl}} &
\ldots &
[\lambda_d+j-d]_{\text{\rm cl}} & \ldots &
[\lambda_d]_{\text{\rm cl}} \\[2mm]
\end{array}
\right|
\]
The $i$th entry on the main diagonal is $[\lambda_i]$ and the index increases
by $1$ per column as we move rightwards in the same row.    The subscript
``cl'' is to remind us that the classes are in ordinary cohomology.    
The theorem below gives an equivariant version of the above formula.

The proof of the equivariant version does not use the ordinary version
of the formula.    In fact, it gives another proof of the ordinary
version by specialization.

Let $u=(u_1,\ldots,u_d)$ in $I(d,n)$.    For $i$, $j$ integers such that
$1\leq i,j\leq d$,  set 
\begin{equation}\label{eudefn}
 u[i,j] = \sum_{k=0}^{j-1}\,\,c(u_i,j,k)\,\,[\lambda_i+j-i-k]	\end{equation}
where $c(u_i,j,k):= (-1)^k\,h_k(\epsilon_{u_i-j+1+k},\ldots,\epsilon_{u_i})$---here
$h_k$ is the ``complete symmetric polynomial,''   the sum of all monomials of
degree $k$ in the elements $\epsilon_{u_i-j+1+k}$, \ldots, $\epsilon_{u_i}$ of 
$H_T^*(\rm{pt})=\bring$.     If $u_i-j+1+k\leq 0$,  then $c(u_i,j,k)$ does not make sense,
but this does not matter since $\lambda_i+j-i-k \geq n-d+1$ and so
$[\lambda_i+j-i-k]=0$ by definition.
\bthm\label{tgiam}\label{tgiambelli}
With notation as above,   given $u=(u_1,\ldots,u_d)$ in $I(d,n)$,  the equivariant
cohomology class $[X(u)]$ is the determinant of the $d\times d$ matrix whose 
$(i,j)\/$th entry is $u[i,j]$.
\ethm
\noindent
This theorem will be deduced in \S\ref{sgiamproof} from
the restriction formula (Theorem~\ref{trestrdet}) and the
injection~(\ref{einjection}).

\section{A determinantal formula for the restriction}\label{srestrdet}
For integers $p$, $k$, $r$,   set
\[\displaystyle{\mu^k_r(p) \coloneq \product_{j=k,\,k+1,\,\ldots,\,r}\epsilon(j,p)}\]
where $\epsilon(j,p)\coloneq \epsilon_j-\epsilon_p$.    This is well-defined 
as an element of the polynomial ring $\bring$
only when $p$, $k$, and $r$ belong to the range $1$, $2$, \ldots, $n$ and $k\leq r$,
but it is convenient to extend the notation somewhat:   if $k=n+1$, the product,
being over an empty index set, is taken to be $1$.
\bthm\label{trestrdet}
Given $u=(u_1,\ldots,u_d)$ and $v=(v_1,\ldots,v_d)$ belonging to $I(d,n)$, 
the restriction $[X(u)]|_v$ of the $T$-equivariant cohomology class $[X(u)]$ of
the Schubert variety $X(u)$ in the Grassmannian $\gdn$
to the $T$-fixed point $e^v$ determined by $v$ equals
\begin{equation}\label{erestrdet}
\displaystyle{ 
\frac{
\left|
\begin{array}{cccccc}
\mu_n^{u_1+1}(v_1) & \mu_n^{u_1+1}(v_2) & \ldots & \mu_n^{u_1+1}(v_j)
	& \ldots & \mu_n^{u_1+1}(v_d) \\[2mm]
\mu_n^{u_2+1}(v_1) & \mu_n^{u_2+1}(v_2) & \ldots & \mu_n^{u_2+1}(v_j)
	& \ldots & \mu_n^{u_2+1}(v_d) \\[2mm]
\vdots & \vdots & \vdots &\vdots & \vdots & \vdots \\[2mm]
\mu_n^{u_i+1}(v_1) & \mu_n^{u_i+1}(v_2) & \ldots &
\mu_n^{u_i+1}(v_j) & \ldots & \mu_n^{u_i+1}(v_d) \\[2mm]
\vdots & \vdots & \vdots &\vdots & \vdots & \vdots \\[2mm]
\mu_n^{u_d+1}(v_1) & \mu_n^{u_d+1}(v_2) & \ldots & \mu_n^{u_d+1}(v_j)
	& \ldots & \mu_n^{u_d+1}(v_d) \\[2mm]
\end{array}
\right| }
{\left|
\begin{array}{cccccc}
1 & 1 & \ldots & 1 & \ldots & 1 \\[2mm]
\epsilon_{v_1} & \epsilon_{v_2} & \ldots & \epsilon_{v_j} & \ldots & \epsilon_{v_d} \\[2mm]
\epsilon_{v_1}^2 & \epsilon_{v_2}^2 & \ldots & \epsilon_{v_j}^2
& \ldots & \epsilon_{v_d}^2 \\[2mm]
\vdots & \vdots & \vdots & \vdots & \vdots & \vdots \\[2mm]
\epsilon_{v_1}^i & \epsilon_{v_2}^i & \ldots & \epsilon_{v_j}^i &
	\ldots & \epsilon_{v_d}^i \\[2mm]
\vdots & \vdots & \vdots & \vdots & \vdots & \vdots \\[2mm]
\epsilon_{v_1}^{d-1} & \epsilon_{v_2}^{d-1} & \ldots & \epsilon_{v_j}^{d-1} &
	\ldots & \epsilon_{v_d}^{d-1} \\[2mm]
\end{array}
\right| }
}\end{equation}
\ethm
\noindent
The denominator in the above expression for $[X(u)]|_v$ is the Vandermonde
determinant which equals
\[ \epsilon(v_2,v_1)\cdot \left(
\epsilon(v_3,v_1)
\epsilon(v_3,v_2)\right) \cdot \ldots \cdot
\left(
\epsilon(v_d,v_1)
\ldots 
\epsilon(v_d,v_{d-1}) \right) \]
The proof of this theorem occupies sections~\ref{srestrproof}
and~\ref{sgrobtorestr}.

\section{Proof of the equivariant Giambelli}\label{sgiamproof}
In this section, Theorem~\ref{tgiambelli} is deduced from Theorem~\ref{trestrdet}.
Because of the injection (\ref{einjection}), it is enough to show that, for an arbitrary $v=(v_1,\ldots,v_d)$ in $I(d,n)$,
the restriction $[X(u)]|_v$ is the determinant of the $d\times d$ matrix 
whose $(i,j){\rm th}$ entry is $u[i,j]|_v$.      We first obtain a 
determinantal formula for $u[i,j]|_v$:   
\begin{equation}\label{edetformula}
u[i,j]|_v = \det(N)/\VDM(v),\end{equation}
where 
\[ \VDM(v):= (\epsilon_{v_2}-\epsilon_{v_1})\cdot
(\epsilon_{v_3}-\epsilon_{v_1})
(\epsilon_{v_3}-\epsilon_{v_2})\cdot\, \ldots\, \cdot
(\epsilon_{v_d}-\epsilon_{v_1})\cdots
(\epsilon_{v_d}-\epsilon_{v_{d-1}})\]
($\VDM$ stands for Vandermonde) and $N$ denotes the following matrix (see \S\ref{trestrdet} for
definition of $\mu_r^k(p)$) 
\[
\left(
\begin{array}{ccccc}
\mu_n^{u_i+1}(v_1)\,(-\epsilon_{v_1})^{j-1} &
\ldots &
\mu_n^{u_i+1}(v_s)\,(-\epsilon_{v_s})^{j-1} &
\ldots &
\mu_n^{u_i+1}(v_d)\,(-\epsilon_{v_d})^{j-1} \\[2mm]
\mu_n^{n-d+3}(v_1) &
%\mu_n^{n-d+3}(v_2) &
\ldots &
\mu_n^{n-d+3}(v_s) &
\ldots &
\mu_n^{n-d+3}(v_d) \\[2mm]
\vdots & %\vdots &
\vdots & \vdots & \vdots & \vdots \\[2mm]
\mu_n^{n-d+r+1}(v_1) &
%\mu_n^{n-d+r+1}(v_2) &
\ldots &
\mu_n^{n-d+r+1}(v_s) &
\ldots &
\mu_n^{n-d+r+1}(v_d) \\[2mm]
\vdots & %\vdots &
\vdots & \vdots & \vdots & \vdots \\[2mm]
\mu_n^{n}(v_1) &
%\mu_n^{n}(v_2) &
\ldots &
\mu_n^{n}(v_s) &
\ldots &
\mu_n^{n}(v_d) \\[2mm]
1 & \ldots & 1 & \ldots & 1 \\[2mm]\end{array} \right)
\]
%\end{equation}
To prove~(\ref{edetformula}),   we substitute for the restrictions of the special classes
on the right side of~(\ref{eudefn}) the determinantal expressions given
by Theorem~\ref{trestrdet}:
\[\begin{array}{rcl}
u[i,j]|_v & =  &
\displaystyle{
\sum_{k=0}^{j-1} c(u_i,j,k) [\lambda_i+j-i-k]|_v } \\[4mm]
	  & =  &
\displaystyle{
\sum_{k=0}^{j-1} c(u_i,j,k)
			\frac{\det(\Num(\lambda_i+j-i-k))}{\VDM(v)} } \\[2mm]
\end{array} \]
where we have written $\Num(\lambda_i+j-i-k)$ for the $d\times d$ matrix whose
determinant is the numerator of the expression for $[\lambda_i+j-i-k]|_v$
given by Theorem~\ref{trestrdet}.     Rows~$2$ through $d$ of $\Num(\lambda_i+j-i-k)$
do not change as $k$ varies:   they are the same as the corresponding ones of the
matrix $N$ in~(\ref{edetformula}).    And the first row of $\Num(\lambda_i+j-i-k)$ is
\[ \left(\mu^{u_i-j+k+2}_n(v_1), \ldots,\mu^{u_i-j+k+2}_n(v_d)\right). \]
So (\ref{edetformula}) follows once we prove
\begin{eqnarray}
\label{eclaimone}
\sum_{k=0}^{j-1}\,\, c(u_i,j,k)\,\,\mu_n^{u_i-j+k+2}(v_s) & = &
			(-\epsilon_{v_s})^{j-1}\mu^{u_i+1}_n(v_s). \end{eqnarray} 
The above identity is the special case $m=j-1$ of the following more
general identity:
for $0\leq m\leq j-1$, 
\begin{eqnarray}
\label{eclaimtwo}
\begin{array}{l}\displaystyle{
\sum_{k=0}^{m}
(-1)^k \,\,\, h_k(\epsilon_{u_i-j+1+k},\ldots,\epsilon_{u_i})
\,\,\mu_n^{u_i-j+k+2}(v_s)\,\, = }\\[2mm]
\displaystyle{
\quad\quad\quad
(-1)^m\,\,\,
h_m(\epsilon_{v_s},\epsilon_{u_i-j+2+m},\ldots,\epsilon_{u_i})\,
\,\,\mu_n^{u_i-j+m+2}(v_s)   }
\end{array}
\end{eqnarray} 
The proof of~(\ref{eclaimtwo}) is by induction on $m$.
First note that it holds for
$m=0$.    Now for the induction step:  assuming that it is true for $m$,
we show it holds for $m+1$.  Taking $(-1)^{m+1} \mu_n^{u_i-j+m+3}(v_s)$
common out of the two terms in the following
\[\begin{array}{l}
(-1)^m
\,\,h_m(\epsilon_{v_s},\epsilon_{u_i-j+2+m},\ldots,\epsilon_{u_i})\,
\,\,\mu_n^{u_i-j+m+2}(v_s)
\,\,\, + \,  \\[2mm]
\quad\quad\quad (-1)^{m+1}
\,\,h_{m+1}(\epsilon_{u_i-j+2+m},\ldots,\epsilon_{u_i})\,
\,\,\mu_n^{u_i-j+m+3}(v_s) \end{array}
\]
we need only show that 
\[\begin{array}{l}
-h_m(\epsilon_{v_s},\epsilon_{u_i-j+2+m},\ldots,\epsilon_{u_i})\cdot
(\epsilon_{u_i-j+m+2}-\epsilon_{v_s}) \\[2mm]
\quad\quad\quad
\,\,+
\,\,h_{m+1}(\epsilon_{u_i-j+2+m},\ldots,\epsilon_{u_i})\\[2mm]
\quad\quad\quad
\quad
\quad
\quad
\,\,=
\,\,h_{m+1}(\epsilon_{v_s},\epsilon_{u_i-j+3+m},\ldots,\epsilon_{u_i}) \end{array}
\]
but this is just the sum of the following two elementary equalities:
\[\begin{array}{rcl}
h_{m+1}(\epsilon_{v_s},\epsilon_{u_i-j+3+m},\ldots,\epsilon_{u_i}) & = & 
h_{m+1}(\epsilon_{v_s},\epsilon_{u_i-j+2+m},\ldots,\epsilon_{u_i}) \\[2mm]
& \quad &-\,\,\epsilon_{u_i-j+2+m}\,\,\,
h_{m}(\epsilon_{v_s},\epsilon_{u_i-j+2+m},\ldots,\epsilon_{u_i})\\[2mm]
h_{m+1}(\epsilon_{v_s},\epsilon_{u_i-j+2+m},\ldots,\epsilon_{u_i}) & = &
h_{m+1}(\epsilon_{u_i-j+2+m},\ldots,\epsilon_{u_i}) \\[2mm]
& \quad &+\,\,\epsilon_{v_s}\,\,\,
h_{m}(\epsilon_{v_s},\epsilon_{u_i-j+2+m},\ldots,\epsilon_{u_i}).\\[2mm]
\end{array}\]

Now that~(\ref{edetformula})
is proved, we proceed with the proof of the theorem.   
If we delete the $1$st row and $s\/$th column of $N$,   the determinant
of the resulting sub-matrix is 
\[ \VDM(v_1,\ldots,\hat{v}_s,\ldots,v_d) :=
\frac{\VDM(v)}{\epsilon({v_s},{v_1})\ldots\epsilon({v_s},v_{s-1})\cdot
\epsilon(v_{s+1},v_s)\ldots\epsilon(v_d,v_s)}.	\]
This follows since the determinant has degree $0+1+\ldots+d-2$ in the
epsilons and is divisible by $\epsilon(v_j,v_i)$ for $1\leq i,j \leq n$,
$i,j\neq s$.
For the determinant of $N$, expanding by the first row, we thus obtain
\[ \det(N) = \sum_{s=1}^d \mu^{u_i+1}_n(v_s)\,\,\epsilon_{v_s}^{j-1}\,\,
		\VDM(v_1,\ldots,\hat{v}_s,\ldots,v_d)	\]
Observe that the right side is the product of the row matrix
\[ \left(\mu^{u_{i}+1}_n(v_1),\ldots,\mu_n^{u_{i}+1}(v_d)\right) \]  with the 
column matrix whose transpose is \[
\left(\epsilon_{v_1}^{j-1}\, \VDM(\hat{v}_1,v_2,\ldots,v_d),
\ldots,
\epsilon_{v_s}^{j-1}\, \VDM(v_1,\ldots,v_{d-1},\hat{v}_d)\right)\]
This means the following for the matrix---let us call it $M$---whose 
$(i,j)$th entry is $u[i,j]|_v$:
$\VDM(v)\,M$  equals
\[
\left(
\begin{array}{ccc}
\mu^{u_{1}+1}_n(v_1) & \ldots & 
\mu^{u_{1}+1}_n(v_d) \\[2mm]
\vdots & \vdots & \vdots \\[2mm]
\mu^{u_{d}+1}_n(v_1) & \ldots & 
\mu^{u_{d}+1}_n(v_d) \\[2mm] \end{array} \right)
\left(
\begin{array}{ccc}
\epsilon^{0}_{v_1}\,\VDM(\hat{v}_1,\ldots,v_d) & \ldots & 
\epsilon^{d-1}_{v_1}\,\VDM(\hat{v}_1,\ldots,v_d) \\[2mm]
\vdots & \vdots & \vdots \\[2mm]
\epsilon^{0}_{v_s}\,\VDM({v}_1,\ldots,\hat{v}_d) & \ldots & 
\epsilon^{d-1}_{v_s}\,\VDM({v}_1,\ldots,\hat{v}_d) \\[2mm]
\end{array}\right)
\]
Since $\product_{s=1}^d \VDM(v_1,\ldots,\hat{v}_s,\ldots,v_d) = \VDM(v)^{d-2}$,
the determinant of the matrix on the right above is $\VDM(v)^{d-1}$.  The matrix
on the left---let us call it $P$---is the numerator in the formula for $[X(u)]|_v$
of Theorem~\ref{trestrdet}.      Taking determinants,  we get
\[ \VDM(v)^d \det(M) = \det(P)\,\VDM(v)^{d-1} \]
%\quad \quad \textrm{
and so by Theorem~\ref{trestrdet}
%} \quad \quad
\begin{flalign*}
\quad\quad\quad\quad\quad
\quad\quad\quad\quad\quad
\det (M) =\det(P)/\VDM(v)= [X(u)]|_v. &&&&\Box
%\quad\quad\quad\quad
%\quad\quad\quad\quad
\end{flalign*}

\section{Proof of the restriction formula}\label{srestrproof}
In this section, Theorem~\ref{trestrdet} is proved. 
Theorem~\ref{tgrobrestr}, which is stated and proved in~\S\ref{sgrobtorestr},
allows us to reduce the proof to combinatorics.   More precisely,
Theorem~\ref{tgrobtorestr} tells us that if we have a Gr\"obner
degeneration of an open piece of the Schubert variety $X(u)$
around the $T$-fixed point $e^v$,   then we can compute
the desired restriction $[X(u)]|_v$.
Such a Gr\"obner degeneration is described in \cite{kr}---indeed
it was the goal of that paper to describe such a degeneration.
We now recall this description.

We identify $\gdn$ as the orbit space for the action on $n\times d$ matrices
of rank~$d$ by the group of
invertible $d\times d$ matrices by multiplication on the right.  
The subset consisting of those matrices in which the submatrix determined by the
rows $v_1,\ldots,v_d$ is the identity matrix gives us an affine $T$-stable
patch of $\gdn$ around the point~$e^v$. 
This patch is an affine space which we denote $\mathbb{A}^v$.
The coordinate function $X(r,j)$ on~$\mathbb{A}^v$
determined by the entry of the matrix in position $(r,j)$, 
$r\not\in v$,  is an eigenvector for~$T$
with character $-(\epsilon_r-\epsilon_{v_j})$.   Thus a
natural way to index these coordinates on~$\mathbb{A}^v$
is by the pairs $(r,c)$, $1\leq r,c\leq n$,
such that $c\in v$ and $r\in\{1,\ldots,n\}\setminus v$---instead
of $X(r,j)$ we write $X(r,v_j)$.   Denote by $\roots^v$ the set of
all such pairs. 

The intersection $Y(u)$
of $X(u)$ with the affine patch~$\mathbb{A}^v$ of $\gdn$ around $e^v$ is
of course a closed subvariety in~$\mathbb{A}^v$.
As proved in \cite[\S5]{kr},  there exist
term orders on the monomials in the variables $X(r,c)$ with
respect to which the initial ideal
of the ideal of functions vanishing on $Y(u)$ is the 
face ideal of a certain simplicial complex~$\srsc$ with vertex set $\roots^v$.
We want to describe the maximal faces of this complex and thereby the complex
itself.

But before we do that,  a digression is necessary.  
In order that specializations to degenerate situations work smoothly,
the correct definition of simplicial complex needs to be adopted.
We do not insist, unlike in \cite[Chapter~II]{stanley:cca} and like in
\cite[Definition~1.4]{millsturm}, on the axiom that
singleton subsets are faces.   More precisely,  here are our definitions:
A {\em simplicial complex\/} is a pair $(V,F)$ of a set $V$
and a set $F$ of subsets of~$V$;  the elements of $V$ are called {\em vertices\/}
and those of $F$ {\em faces\/};  the following axioms hold:
(1)~the empty subset of $V$ is a face, and (2)~a subset of a face is a face.
Because of (2) we may replace~(1) by (1'): $F$ is non-empty.

Given a simplicial complex $(V,F)$, its {\em (Stanley-Reisner) face ring\/} $\srfr$ is
defined as follows:  consider
the polynomial ring, over some implicit base, in a set of variables indexed 
by~$V$---we abuse notation and let $V$ itself denote the set of variables;  the
linear span of monomials whose support is not contained in any face
forms an ideal---let us call it the {\em face ideal\/} (or should it be
the {\em non-face ideal\/}?);
the quotient of the polynomial ring by the face ideal is $\srfr$.
It is readily seen that the face ideal is the intersection, 
over all maximal faces,
of the ideal generated by the variables in the complement of that face.

The digression being over, we now start on the description of the
simplicial complex~$\srsc$.
Denote by $\pos^v$ the subset of $\roots^v$
consisting of those pairs $(r,c)$ for which $r>c$.
The element $u$ of $\idn$ determines as follows a subset $\mon_u^v$ of $\pos^v$ 
with the following property:
writing $\mon_u^v =\{(r_1,c_1),\ldots,(r_k,c_k)\}$,  we have
$u=\left(v\setminus\{c_1,\ldots,c_k\}\right)\cup\{r_1,\ldots,r_k\}$.
To define $\mon_u^v$, proceed by induction on $d$.
Let $i$, $1\leq i\leq d$, be the largest such that $v_i\leq u_1$. 
Set $v'=v\setminus\{v_i\}$ and $u'=u\setminus\{u_1\}$.   Then $v'\leq u'$
and $\mon_{u'}^{v'}$ is defined by induction.   Set
\[ \mon_u^v=\left\{\begin{array}{ll}
 	\mon_{u'}^{v'}\cup\{(u_1,v_i)\} & \text{if $u_1\neq v_i$}\\ 
 	\mon_{u'}^{v'}			& \text{if $u_1=v_i$}\\
	\end{array}\right. \] 

\begin{figure}%\label{feone}
\begin{center}
%\center{\mbox{\epsfbox{feoneold.eps}}} 
\mbox{\epsfig{file=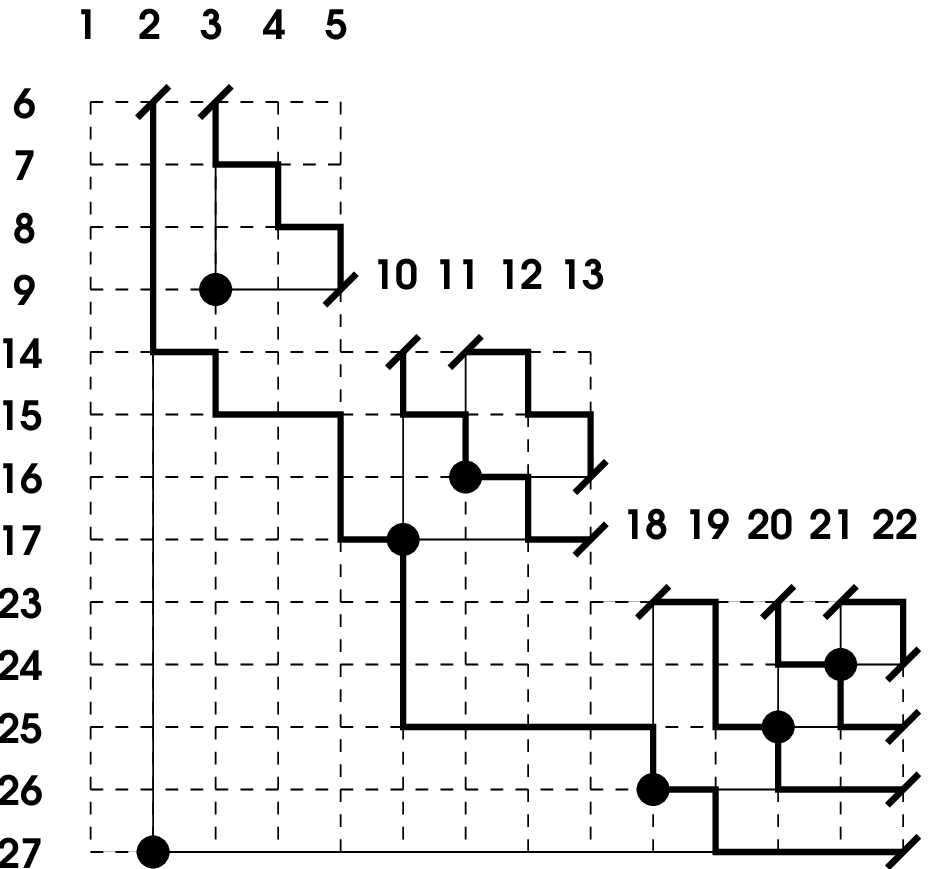,height=5.2cm,width=6cm}}
\caption{\label{feone}A tuple of non-intersecting lattice paths as in 
	Example~\ref{eone}}
\end{center}
\end{figure}

We draw---see Example~\ref{eone} and
Figure~\ref{feone} below---a grid with the
elements of $\pos^v$ being the lattice points---in the notation $(r,c)$,
the $r$ is suggestive of row index and $c$ of column index.    The solid
dots in the figure denote the points of $\mon_u^v$.  From each solid
dot $\beta$ we draw a vertical line and a horizontal line.
Let $\bstart$ and $\bfinish$ denote respectively the points
where the vertical and the horizontal lines meet the boundary.
In Figure~\ref{feone} for example $\bstart=(14,11)$ and
$\bfinish=(16,13)$ for $\beta=(16,11)$.

A {\em lattice path\/} between a pair of such points $\bstart$
and $\bfinish$ is a sequence $\alpha_1,\ldots,\alpha_q$
of elements of $\pos^v$ with $\alpha_1=\bstart$,
$\alpha_q=\bfinish$,  and for $j$, $1\leq j\leq q-1$,
writing $\alpha_j=(r,c)$,  $\alpha_{j+1}$ is either
$(r',c)$ or $(r,c')$, where $r'$ (respectively $c'$) is the smallest
integer not in $v$ (respectively in $v$)
and greater than $r$ (respectively $c$).     If  $\bstart=(r,c)$ and
$\bfinish=(R,C)$,   then $q=(R-r)+(C-c)+1$.

Let us write $\mon_u^v=\{\beta_1,\ldots,\beta_p\}$.  Consider
the set of all $p$-tuples of paths $\path=(\path_1,\ldots,\path_p)$,
where $\path_j$ is a lattice path between $\bstartj$
and $\bfinishj$,  and no two $\path_j$ intersect.
A particular such $p$-tuple is shown in Figure~\ref{feone}.
Such $p$-tuples form an indexing set
for the maximal faces of the
simplicial complex $\srsc$:   
to the $p$-tuple $\path=(\path_1,\ldots,\path_p)$ 
corresponds
the maximal face 
$\path_1\cup\cdots\cup\path_p\cup
\left(\roots^v\setminus\pos^v\right)$: %\footnote{%See footnote~\ref{fdegenerate}
%for the meaning of this bijection
in the degenerate case when $\mon^v_u$ 
is empty (which happens only if $u=v$), there is a unique maximal
face, namely $\roots^v\setminus\pos^v$.%}
\bexample\label{eone}
Let $d=14$, $n=27$,  
\begin{eqnarray*}
v&=&(1,2,3,4,5,10,11,12,13,18,19,20,21,22),
\mbox{\begin{rm} and \end{rm}}\\
u&=&(1,4,5,9,12,13,16,17,19,22,24,25,26,27), 
\mbox{\begin{rm} so that \end{rm}}\\
\mon_u^v &=& \left\{(9,3), (16,11), (17,10), (24,21), (25,20),  (26,18),
			(27,2) \right\}.
\end{eqnarray*}
Figure~\ref{feone} shows a particular tuple of non-intersecting
lattice paths.\eexample

Consider now the subvariety of $\mathbb{A}^{v}$
defined by the face ideal
of the
complex~$\srsc$.
It is a union of coordinate planes.   There is one plane
for each maximal face and it is defined by the
vanishing of the coordinates corresponding to the vertices
in the complement of that face.     
For a maximal face $\face$ corresponding to $(\path_1,\ldots,\path_p)$, 
denote by $m_\face$ the product,
over all $(r,c)$
in $\pos^v\setminus\left(
\path_1\cup\cdots\cup\path_p\right)$,
of $\epsilon(r,c):=\epsilon_r-\epsilon_c$.
It follows from Theorem~\ref{tgrobrestr}
that the
restriction $[X_u]|_v$ is the sum $\sum m_\face$
as $\face$ varies over all maximal
faces. %\footnote{\label{fdegenerate}
The last assertion holds also
in the degenerate case $u=v$:
then $\mon^v_u$ is empty;  $\srsc$ has only one maximal face,
namely $\roots^v\setminus \pos^v$; and $\sum m_\face$ is the product 
over $(r,c)$ in $\pos^v$ of $(\epsilon_r-\epsilon_c)$.   In particular,
\begin{itemize}
\item
If $u=v=(1,2,\ldots,d)$,  then $\srsc$ has only the empty face,
and $\sum m_\face$ is the product 
over $(r,c)$ in $\pos^v=\roots^v$ of $(\epsilon_r-\epsilon_c)$.
\item
If $u=v=(n-d+1,\ldots,n)$,  then the unique maximal face of $\srsc$
is $\face=\roots^v$
and $m_\face$, being the product over an empty index set, equals $1$.
\end{itemize}%}

\begin{figure}[h]%\label{fetwo}
\begin{center}
\mbox{\epsfig{file=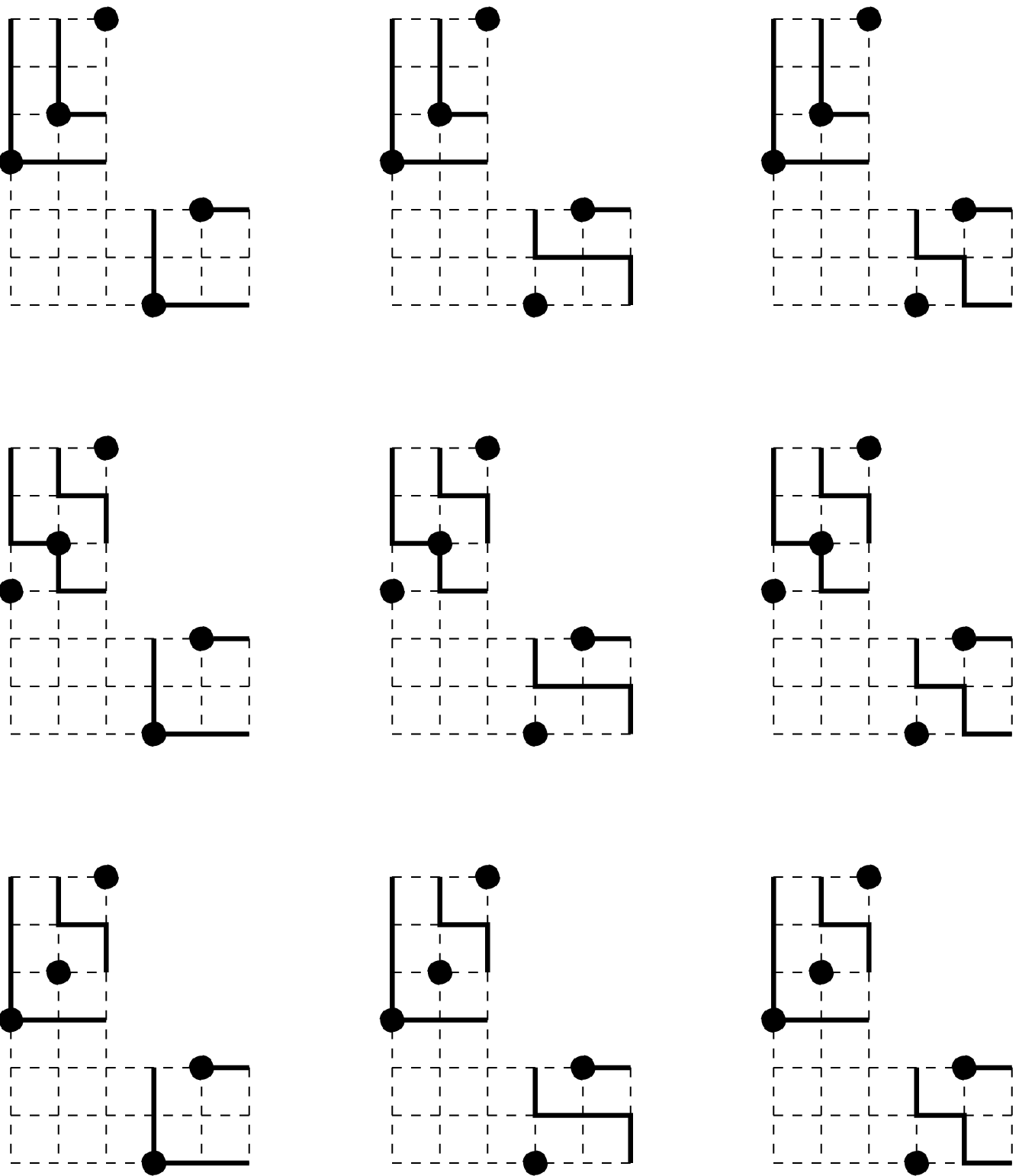,height=6cm,width=6cm}}
\caption{\label{fetwo}All the tuples of non-intersecting lattice paths as in 
	Example~\ref{etwo}}
\end{center}
\end{figure}

\bexampletwo\label{etwo}
This example is simple enough so we can easily draw all
possible tuples of non-intersecting lattice paths.
Let $d=6$, $n=13$,
\[\begin{array}{c}
v=(1,2,3,8,9,10), \text{ and }
u=(4,6,7,10,11,13).\\  
\text{Then }\mon_u^v=\{(4,3), (6,2),
(7,1), (11,9), (13,8)\}.\end{array}\]
Figure~\ref{fetwo} shows all 
the $5$-tuples of 
non-intersecting lattice paths---there are~$9$ of them.  
Writing $\epsilon(r,c)$ for $\epsilon_r-\epsilon_c$, 
\[\begin{array}{rcl}
[X_u]|_v& = &
\epsilon(11,1)
\epsilon(11,2)
\epsilon(11,3)
\epsilon(12,1)
\epsilon(12,2)
\epsilon(12,3)
\epsilon(13,1)
\epsilon(13,2)
\epsilon(13,3)\cdot\\
%\quad\quad\quad
&&
\left[
\epsilon(12,9)\epsilon(12,10)+\epsilon(13,8)\epsilon(12,10)+\epsilon(13,8)\epsilon(13,9)
\right]\cdot\\
%\quad\quad\quad\quad\quad\quad
&&\left[\epsilon(5,3)+\epsilon(6,2)+\epsilon(7,1)\right].\Box
\end{array}
\]\eexampletwo

Thus the proof of the restriction formula is reduced to the combinatorial
problem of establishing 
\begin{equation}\label{ecomb}
\sum m_\face=E(u,v)\end{equation}
where $E(u,v)$ stands for the expression~(3).
Whether this problem admits of an elegant
solution by means of a Lindstrom-Gessel-Viennot type argument we do not know.
What follows is an elementary argument based on induction.

Proceed by induction on $d$.   The case $d=1$ being easily verified,
let $d\geq2$.    The strategy of the proof is this.  In the first part,
we work with $\sum m_\face$ and express it in terms of ``smaller''
$\sum m_\face$---those attached to simplicial complexes $\srscprime$ for
elements $u'\geq v'$ in~$I(d-1,n)$.     By the induction hypothesis,
equality (\ref{ecomb}) applies to these
smaller $\sum m_\face$, so that  we get an expression for
$\sum m_\face$ in terms of $E(u',v')$---the precise expression
is in~(\ref{eend}) below.
In the second part,  we will algebraically manipulate the expression
(\ref{erestrdet}) for $E(u,v)$ to express it in terms
of~$E(u',v')$.   
The resulting expression will turn out to be the same as that for
$\sum m_\face$ obtained in the first part.
This will finish the proof.

\begin{figure}%\label{feone}
\psfrag{ldone}{$(\ell_d,v_{d-1})$}
\psfrag{ldtwo}{$(\ell_d,v_{d})$}
\psfrag{lrone}{$(\ell_{r+1},v_{d-1})$}
\psfrag{lrtwo}{$(\ell_{r+1},v_{d})$}
\psfrag{vone}{$v_{d-1}$}
\psfrag{vtwo}{$v_{d}$}
\psfrag{ur}{$(u_{r},v_{d})=\beta_{p-d+r}$}
\psfrag{urtwo}{$=\beta_{p-d+r}(\textrm{finish})$}
\psfrag{ud}{$(u_{d},v_{d})=\beta_{p}(\textrm{finish})$}
\psfrag{urone}{$(u_{r+1},v_{d})$}
\psfrag{vdots}{$\vdots$}
\psfrag{uronetwo}{$=\beta_{p-d+r+1}(\textrm{finish})$}
\psfrag{udone}{$(u_{d-1},v_{d})=\beta_{p-1}(\textrm{finish})$}
\psfrag{ltwo}{$(\ell_{d-1},v_{d-1})$}
\psfrag{lone}{$(\ell_{d-1},v_{d})$}
\begin{center}
\mbox{\epsfig{file=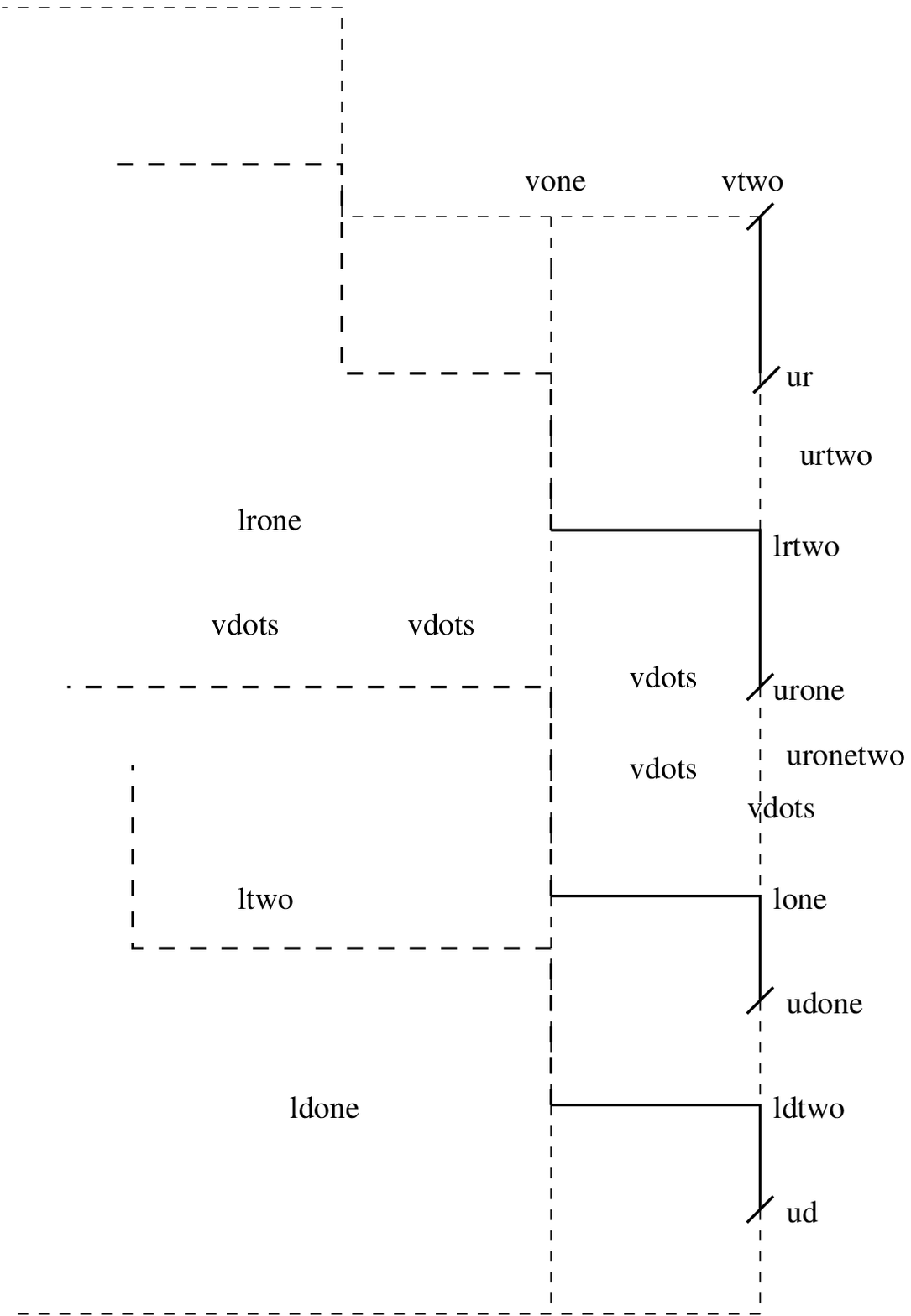,height=8.0cm,width=5.8cm}}
\caption{\label{fproof} Partitioning the maximal faces of $\srsc$ into $S_{\seq}$}

\end{center}
\end{figure}

So first consider $\sum m_\face$.
Let $r$ be the integer, $1\leq r\leq d$, such that
$u_{r-1}<v_d\leq u_r$.    
Write as
before $\mon^v_u=\{\beta_1,\ldots,\beta_p\}$.    It is easy to see
that
\[\begin{array}{c}
\beta_p(\textrm{finish})=(u_d,v_d), \ 
\beta_{p-1}(\textrm{finish})=(u_{d-1},v_d), \ \ldots, \\[2mm]
\quad\quad\quad\beta_{p-d+r+1}(\textrm{finish})=(u_{r+1},v_d);\\[2mm]
\textrm{furthermore,
$\beta_{p-d+r}(\textrm{finish})=(u_{r},v_d)$ unless $u_r=v_d$.}\end{array}  \]
Figure~\ref{fproof} depicts the situation.

Partition the set $S$ of the $p$-tuples of paths $(\path_1,\ldots,\path_p)$
(those indexing the maximal faces of $\srsc$)
into subsets $S_{\seq}$ indexed by sequences $\seq=(\ell_{r+1},\ldots,\ell_d)$
of integers such that $u_r<\ell_{r+1}\leq u_{r+1}$,  $u_{r+1}<\ell_{r+2}\leq
u_{r+2}$, \ldots, $u_{d-1}<\ell_d\leq u_d$:       the subset $S_{\seq}$ consists
of all those in which the segment joining $(\ell_j,v_{d-1})$ and 
$(\ell_j,v_d)$ is part of the path $\path_{p-d+j}$ for every $j$,
$r+1\leq j\leq d$. 
That the $S_{\seq}$ form a partition of~$S$ is readily seen.
Letting $S_{\seq}$ also denote
the corresponding partition of the maximal faces of $\srsc$, we obtain
\begin{equation}\label{erestrproofone}
\sum m_\face = \sum\limits_{\seq} \sum\limits_{f\in S_{\seq}} m_\face
\end{equation}

For the moment, let us fix a certain $\seq$.  Set 
$u':=(u_1,\ldots,u_{r-1},l_{r+1},\ldots,l_d)$  and
$v':=(v_1,\ldots,v_{d-1})$ (although $u'$ depends on $\seq$,
we still write only $u'$ and not $u'(\seq)$).
Then $u'\geq v'$. 
We want to use the induction hypothesis to express $\sum_{\face\in S_{\seq}} m_\face$ in terms of $E(u',v')$.      Towards this, we make two observations.
First,  the factor
\[\mu_n^{u_d+1}(v_d)\cdot
\product_{r\leq i<d}\mu_{\ell_{i+1}-1}^{u_i+1}(v_d)\]
is common to all the terms in $\sum_{\face\in S_{\seq}}m_\face$.   Second,  the
integer $v_d$ does not occur as a row or column index if we restrict attention
to the first $d-1$ columns of Figure~\ref{fproof} (which tells us that
$E(u',v')$ needs to be adjusted to take care of this).

Consider $E(u',v')$---this is the expression~(3) with $u$ and $v$
replaced respectively by $u'$ and $v'$.
In the matrix whose determinant is the numerator of~$E(u',v')$, 
the entry in position $(i,j)$ where $i\leq r-1$
has $\epsilon(v_d,v_j)$ %=(\epsilon_{v_d}-\epsilon_{v_j})$
occurring as a
factor---this factor does not 
occur if $i\geq r$.     Denote by $E(u',v';v_d)$ the modified
expression where %in~(\ref{erestrdet})
the factors 
$\epsilon({v_d},{v_j})$ are taken out---more precisely,
$E(u',v';v_d):=$
\newcommand{\dsfrac}{\ds\frac}
{\begin{equation}\label{erestrdetmod}
\dsfrac{
\left|
\begin{array}{ccccc}
\dsfrac{\mu_n^{u_1+1}(v_1)}{\epsilon(v_d,v_1)}
& \ldots & \dsfrac{\mu_n^{u_1+1}(v_j)}{\epsilon(v_d,v_j)}
	& \ldots & \dsfrac{\mu_n^{u_1+1}(v_{d-1})}{\epsilon(v_d,v_{d-1})}\\[2mm]
\vdots & \vdots &\vdots & \vdots & \vdots \\[2mm]
\dsfrac{\mu_n^{u_{r-1}+1}(v_1)}{\epsilon(v_d,v_1)}
& \ldots & \dsfrac{\mu_n^{u_{r-1}+1}(v_j)}{\epsilon(v_d,v_j)}
	& \ldots & \dsfrac{\mu_n^{u_{r-1}+1}(v_{d-1})}{\epsilon(v_d,v_{d-1})}\\[4mm]
{\mu_n^{l_{r+1}+1}(v_1)}
& \ldots &{\mu_n^{l_{r+1}+1}(v_j)}
	& \ldots & {\mu_n^{l_{r+1}+1}(v_{d-1})}\\[2mm]
\vdots & \vdots &\vdots & \vdots & \vdots \\[2mm]
{\mu_n^{l_{d}+1}(v_1)}
& \ldots &{\mu_n^{l_{d}+1}(v_j)}
	& \ldots & {\mu_n^{l_{d}+1}(v_{d-1})}\\[2mm]
\end{array}
\right| }
{
\epsilon(v_2,v_1)\cdot \left(
\epsilon(v_3,v_1)
\epsilon(v_3,v_2)\right) \cdot \ldots \cdot
\left(
\epsilon(v_{d-1},v_1)
\ldots 
\epsilon(v_{d-1},v_{d-2}) \right)
}\end{equation}}
%%%%%%%%%%%%%%%
Looking at Figure~\ref{fproof} and using the induction hypothesis, we get
{\begin{equation}\label{eend}
\sum\limits_{\face\in S_{\seq}} m_\face \, = \,
E(u',v';v_{d})\cdot
\mu^{u_{d}+1}_{n}(v_d)\cdot\!
\displaystyle\prod_{r\leq i<d}
\mu_{\ell_{i+1}-1}^{u_i+1}(v_d).
\end{equation}}

We are done with the first half of the proof.   Namely, we are finished
with the ``combinatorial side'' $\sum m_\face$ of Equation~(\ref{ecomb}).  
Next we turn to the ``algebraic side'' $E(u,v)$ and show
that it too is a sum of terms indexed by the sequences~$\seq$.  
We show that the term
corresponding to a sequence $\seq$ equals the
right hand side of~(\ref{eend}).   This clearly suffices to complete the proof.

By definition, $E(u,v)$ is the expression (\ref{erestrdet}).   
The entries in the last
column of the numerator vanish in rows $i$ for $1\leq i \leq r-1$ because
$u_i+1\leq u_{r-1}+1\leq v_d$ and $0=\epsilon_{v_d}-\epsilon_{v_d}$ occurs as a factor in
$\mu_n^{u_i+1}(v_d)$.     We would like to kill also the entries in
rows $r$ through $d-1$.   To this end,  subtract from row $i$,
$r\leq i\leq d-1$,
the multiple of row~$d$ by $\mu_{u_d}^{u_i+1}(v_d)$.   
The entry in position $(i,j)$ %for $r\leq i\leq d-1$
then becomes
\[
\mu_n^{u_i+1}(v_j)-\mu_{u_d}^{u_i+1}(v_d)\mu_n^{u_d+1}(v_j)
=\mu_n^{u_d+1}(v_j)\left(\mu_{u_d}^{u_i+1}(v_j)-\mu_{u_d}^{u_i+1}(v_d)\right).\]
In particular, all the entries in the last column except the one on
row $d$ are $0$.

The factor $\mu_n^{u_d+1}(v_j)$ being now
common to all the entries in column $j$,
let us take these factors out of every column. 
The resulting entries in the last column are all zero except the one on row $d$
which is $1$.
The numerator therefore reduces to the determinant of the 
submatrix of the first $d-1$ rows and columns.

Let us take the
factors $(\epsilon_{v_d}-\epsilon_{v_1})$, \ldots, $(\epsilon_{v_d}-\epsilon_{v_{d-1}})$
out of the denominator and distribute
them thus:  divide by $(\epsilon_{v_d}-\epsilon_{v_j})$ every entry in column
$j$ of the determinant in the numerator.  
After these manipulations, expression
(\ref{erestrdet}) looks like this:
\begin{equation}
\label{eintermediate}
\prod\limits_{j=1}^d\mu_n^{u_d+1}(v_j)
\frac{\det(A)}{
%[(v_2-v_1)]\,[(v_3-v_2)(v_3-v_1)]
%\cdots[(v_{d-1}-v_1)\cdots(v_{d-1}-v_{d-2})]}
\epsilon(v_2,v_1)%\cdot
\left(
\epsilon(v_3,v_1)
\epsilon(v_3,v_2)\right)
%\cdot
\cdots
%\cdot
\left(
\epsilon(v_{d-1},v_1)
\ldots 
\epsilon(v_{d-1},v_{d-2}) \right)
}\end{equation}
where $A:=(A_{ij})$ is the $d-1\times d-1$ matrix whose entry at position $(i,j)$ is
\begin{equation}
\label{eentry}
A_{ij} =\frac{\mu_{u_d}^{u_i+1}(v_j)-\mu_{u_d}^{u_i+1}(v_d)}
	{\epsilon_{v_d}-\epsilon_{v_j}}.	\end{equation}

Now apply the following row operations to the matrix $A$: subtract
row~$r+1$ from row $r$, \ldots, row $d-1$ from row $d-2$. 
To get a handle on the resulting matrix---let us call it $B$---we
use the following equation
which is proved readily by induction:
for positive integers $a\leq e\leq b$, $c$, and $f$, we have
\begin{equation*}%\label{ename}
\begin{array}{rcl}
\dsfrac{\mu_b^a(c)-\mu_b^a(f)}{\epsilon_f-\epsilon_c} & = &
\mu_b^{a+1}(c) + \mu_b^{a+2}(c)\mu^a_a(f)+\cdots+\mu^e_b(c)\mu_{e-2}^a(f)+\\[2mm]
& & \quad\quad\quad\quad\quad\quad\quad\quad\quad\quad\quad\quad\quad\quad\quad
\dsfrac{\mu^e_b(c)-\mu^e_b(f)}{\epsilon_f-\epsilon_c} \\[2mm]
& = &
\left(\sum\limits_{a-1< \ell \leq e-1}\mu^{\ell+1}_b(c)\cdot\mu^{a}_{\ell-1}(f)\right)
+\dsfrac{\mu^e_b(c)-\mu^e_b(f)}{\epsilon_f-\epsilon_c}
\end{array}
\end{equation*}
Applying this to (\ref{eentry}) for $i$ such that $r\leq i\leq d-1$ and
$e=u_{i+1}+1$,  we get   %The last term on the right of (\ref{ename}) is then
\begin{equation*}\begin{array}{rcl}
%\label{eentrytwo}
\dsfrac{\mu_{u_d}^{u_{i}+1}(v_j)-\mu_{u_d}^{u_{i}+1}({v_d})}{\epsilon_{v_d}-\epsilon_{v_j}}  & = &  
\left(\sum\limits_{u_{i}< \ell_{i+1} \leq {u_{i+1}}}\mu^{\ell_{i+1}+1}_{u_d}(v_j)\cdot\mu^{u_i+1}_{\ell_{i+1}-1}({v_d})\right)+\\[5mm]
&&%\dsfrac{\mu^e_{u_d}(v_j)-\mu^e_{u_d}(d)}{\epsilon_d-\epsilon_{v_j}}
	\quad\quad\quad\quad\dsfrac{\mu_{u_d}^{u_{i+1}+1}(v_j)-\mu_{u_d}^{u_{i+1}+1}(v_d)}
	{\epsilon_{v_d}-\epsilon_{v_j}}\end{array}\end{equation*}
where we have written $\ell_{i+1}$ rather than just $\ell$ for the running
index in the sum.   Note that the second term on the right is 
precisely the entry at position $(i+1,j)$ of $A$ for $r\leq i\leq d-2$ and vanishes for $i=d-1$.     Thus the entry at position $(i,j)$ of the matrix $B$ looks
like this:
\begin{equation*}
B_{ij}=\left\{
\begin{array}{ll}
\dsfrac{\mu_{u_d}^{u_i+1}(v_j)}{\epsilon(v_d,v_j)}
&\quad\textrm{if $i<r$}\\[4mm]
\sum\limits_{u_{i}<\ell_{i+1}\leq u_{i+1}}{\mu_{u_d}^{l_{i+1}+1}(v_j)}\cdot
{\mu_{\ell_{i+1}-1}^{u_{i}+1}(v_d)}
&\quad\textrm{if $i\geq r$.}\end{array}\right.\end{equation*}
By the multilinearity of the determinant,
we see that $\det(B)$ (which equals $\det(A)$, since
$B$ was obtained from $A$ by elementary
row operations) equals the sum over 
${\seq=(\ell_{r+1},\ldots,\ell_d)}$ of 
\begin{equation}\label{elast}
%\sum\limits_{\seq=(\ell_{r+1},\ldots,\ell_d}
%{\mu_{u_d}^{l_{i+1}+1}(v_j)}\cdot
\product\limits_{r\leq i<d}
{\mu_{\ell_{i+1}-1}^{u_{i}+1}(v_d)}
\cdot\left|
\begin{array}{ccccc}
\dsfrac{\mu_{u_d}^{u_1+1}(v_1)}{\epsilon(v_d,v_1)}
& \ldots & \dsfrac{\mu_{u_d}^{u_1+1}(v_j)}{\epsilon(v_d,v_j)}
	& \ldots & \dsfrac{\mu_{u_d}^{u_1+1}(v_{d-1})}{\epsilon(v_d,v_{d-1})}\\[2mm]
\vdots & \vdots &\vdots & \vdots & \vdots \\[2mm]
\dsfrac{\mu_{u_d}^{u_{r-1}+1}(v_1)}{\epsilon(v_d,v_1)}
& \ldots & \dsfrac{\mu_{u_d}^{u_{r-1}+1}(v_j)}{\epsilon(v_d,v_j)}
	& \ldots & \dsfrac{\mu_{u_d}^{u_{r-1}+1}(v_{d-1})}{\epsilon(v_d,v_{d-1})}\\[4mm]
{\mu_{u_d}^{l_{r+1}+1}(v_1)}
& \ldots &{\mu_{u_d}^{l_{r+1}+1}(v_j)}
	& \ldots & {\mu_{u_d}^{l_{r+1}+1}(v_{d-1})}\\[2mm]
\vdots & \vdots &\vdots & \vdots & \vdots \\[2mm]
{\mu_{u_d}^{l_{d}+1}(v_1)}
& \ldots &{\mu_{u_d}^{l_{d}+1}(v_j)}
	& \ldots & {\mu_{u_d}^{l_{d}+1}(v_{d-1})}\\[2mm]
\end{array}\right|.\end{equation}
Substitute this into equation~(\ref{eintermediate}).   Multiplying the factor $\mu_n^{u_{d}+1}(v_j)$, for $1\leq j\leq d-1$, in equation~(\ref{eintermediate}) into all the entries in column $j$ of the determinant in (\ref{elast}) yields the determinant in the numerator of $E(u',v';v_d)$ in (\ref{erestrdetmod}).   It should now be clear
that $E(u,v)$ is the sum over $\seq$ of the
right side of (\ref{eend}),   and the
proof of the restriction formula (Theorem~\ref{trestrdet}) is finally over.

\section{Gr\"obner degeneration computes restriction}\label{sgrobtorestr}\label{sgrobrestr}
\renewcommand{\tilde}{\widetilde}
The goal of this section is to state
and prove Theorem~\ref{tgrobrestr} below,
which was used 
in the proof in~\S\ref{srestrproof} of the restriction formula
(Theorem~\ref{trestrdet}).     
%We thank the referee for the proof given below of Theorem~\ref{tgrobrestr}.
As pointed out to us by the referee,  Theorem~\ref{tgrobrestr}
is well known and can be deduced from results in the literature.

%But nevertheless we include a proof the theorem here for the sake
%of completeness.  
%In content and spirit, this section is different from the rest of the
%paper.  

The assumptions and notations of
the previous sections are annulled now. 
Fix a torus $T:=(\mathbb{C}^*)^m$.   Let $Z$ be a non-singular
complex projective variety of dimension~$d$ on which there is an
algebraic action
of $T$ with finitely many fixed points.  Then,
by Bialynicki-Birula~\cite{bb},
$Z$ admits an equivariant algebraic
cell decomposition,  and
around each $T$-fixed point there is a $T$-stable open subset~$U$
of $Z$ that is isomorphic to a $T$-module (the fixed point
of course corresponds to $0$ in the module).

Let $Y$ be a $T$-stable irreducible subvariety of $Z$, $y$ a $T$-fixed
point on $Y$,  and $U\simeq \mathbb{C}^d$
a $T$-stable open subset of $Z$ containing~$y$ as above.
Let $X_1,\ldots,X_d$
be coordinates on~$U$ that are eigenvectors for $T$---since $y$ is 
an isolated fixed point,  no coordinate has trivial character.
Choose some term order on the set of monomials in the coordinates,
and suppose that $J$ is the initial ideal,
with respect to this order, of the ideal of functions on~$U$ vanishing
on~$Y\cap U$.

\newcommand{\idealp}{{\mathfrak p}}
Since $J$ is a monomial ideal, it has a primary decomposition
consisting of monomial ideals.    Let 
$\cap_{i=1}^p J_i$ be the intersection of the minimal primary
components (we are throwing away the embedded components).
The radical $\idealp_i$ of $J_i$ is of the form 
$(X_1^{a_{1i}},\ldots,X_d^{a_{di}})$ where each $a_{ji}$ is either~$0$
or~$1$ 
(exactly $\dim Y$ of the $a_{ji}$ equal $0$ for each $i$).
The scheme $M$ defined by $\cap_{i=1}^p J_i$ is the union $\cup_{i=1}^p M_i$
of the schemes~$M_i$ defined by~$J_i$.   Let $\ell_i$ be the length
of $R_{\idealp_i}/JR_{\idealp_i}=R_{\idealp_i}/J_iR_{\idealp_i}$
 where $R:=\mathbb{C}[X_1,\ldots, X_d]$.
\bthm\label{tgrobtorestr}\label{tgrobrestr} With hypothesis and notation as above,
the restriction to the fixed point~$y$
of the equivariant
cohomology class $[Y]$ of the subvariety $Y$ in the equivariant integral
cohomology ring of $Z$ is given by 
\begin{equation}\label{etgrob}
[Y]|_y = \sum\limits_{i=1}^p\ell_i\product\limits_{j=1}^{d}
\chi_j^{a_{ji}}
\end{equation}
where $-\chi_1$, \ldots, $-\chi_d$ are respectively the characters of $X_1$,\ldots, $X_d$.\ethm
\noindent
\bproof
We thank the referee for indicating how the theorem can be deduced
as follows from known results.    The restriction of the equivariant
Chow class of $Y$ to the fixed point $y$ factors through the
restriction to the open set $U$.     The fact that the class of $Y\cap U$
and of $M$ in the equivariant Chow ring of $U$ are the same and
equal to the right hand side of~\eqref{etgrob} can be found in 
%The equivariant chow class 
%In the notation to be established below,  that the equivariant class
%of $Y\cap U$ in the torus-representation $U$ is the same as the
%right hand side of the equation in Theorem~\ref{tgrobrestr} is 
%the rational equivalence of the classes of $E_n(Y\cap U)$ and $E_n(M)$
%in $E_n(U)$ can be found in
any number of references under the
heading of ``equivariant Chow'', or ``multidegree'', or ``equivariant
Hilbert polynomials'', or ``equivariant multiplicity.''   See, for
example,  \cite[Notes to Chapter~8]{millsturm};   the fact that the
``multidegrees''
in the above reference are equivariant cohomology classes is asserted in
Proposition~1.19 of~\cite{kms},  where multidegrees are identified as
being equivariant Chow classes.\eproof
\ignore{
\bproof
We start by recalling some general facts about equivariant cohomology.
The equivariant cohomology $H_T^*(X)$ of a $T$-space $X$ is defined as the
cohomology of $E(X):=E\times_T X$ where $E\to B$ is a universal $T$-bundle.
Since $E$ and $B$ are infinite dimensional
objects,  we do not work with them directly,
but instead with
finite dimensional approximations to them.  We need to ensure that
the approximations are big enough so that no relevant information is lost,
and this is the aim of the next three paragraphs.

For $n$ a positive integer, set $\bE_n:=(\bc^{n+1}\setminus 0)^m$ 
and $\bb_n:=(\bp^n)^m$.  Then 
$\bE_n\lr \bb_n$ is a $T$-bundle, the $m$-fold product of the 
tautological $\mathbb{C}^*$-bundle over $\bp^n$. The inclusion of $\bc^{n+1}$ 
into $\bc^{n+2}=\bc^{n+1}\oplus \bc$ induces inclusions
$\bb_n \subset \bb_{n+1}$ and 
$\bE_n \subset \bE_{n+1}$.
The bundle projection 
$\bE_{n+1}\lr \bb_{n+1}$ 
restricts to 
$\bE_n\lr\bb_n$.   
Therefore, setting $E:=\cup_{n\geq 1} \bE_n$ and
$B:=\cup_{n\geq 1}\bb_n=(\bp^\infty)^m$,
we get a $T$-bundle projection $E\to B$. Since $E$ is contractible,
it follows that $\pi:E\lr B$
is a universal $T$-bundle.
Thus $B$ is a universal 
classifying space for $T$-bundles.

For a $T$-space $X$,  set
$\bE_n(X):=\bE_n\times_T X$.  We have the associated bundle
$\bE_n(X)\to\nobreak \bb_n$ with fiber $X$. 
Since $\bE_n$ contains the $(2n+1)$-skeleton of $E$,  it follows that the 
inclusion $\bE_n(X)\hookrightarrow E(X)$ induces isomorphisms in equivariant cohomology  
$H^k_T(X)\lr H^k(\bE_n(X))$ in degrees $k$, $1\leq k\leq 2n$.   

Fix an integer $n$ larger than twice the (complex) codimension of $Y$ in $Z$.
We will use $E_n(\cdot)$ instead of $E(\cdot)$ to compute the equivariant
cohomology class of $Y$.   
Since the two cohomologies are isomorphic in all degrees in
which the relevant classes live,    this is justified.

Since $E_n(\{y\})\simeq B_n=(\mathbb{P}^n)^m$,  we may identify the singular
and Chow cohomology rings of $E_n(\{y\})$.    Further, since $E_n(U)$ is a
vector bundle over $E_n(\{y\})$,   their Chow rings are naturally isomorphic.
We will show that the classes $[E_n(Y\cap U)]$
and $[E_n(M)]$ in the Chow ring of $E_n(U)$ are equal. 
Since $[E_n(Y\cap U)]$---or, more precisely, its image in the cohomology ring
of $E_n(\{y\})$---is the restriction on the left side of (\ref{etgrob}),
it then remains only to show that $[E_n(M)]$ equals the right side
of~(\ref{etgrob}).

Make $U\times\mathbb{A}^1$ into a $T$-variety with
$T$ acting diagonally on the product and trivially on $\mathbb{A}^1$.
There is a certain $T$-stable subscheme $\tilde{V}$ of $U\times\mathbb{A}^1$
such that the composition $\tilde{V}\to\nobreak\mathbb{A}^1$ is a flat morphism with
the fiber over $0$ being~$M$ and fiber over any other point being~$Y\cap U$---see,
for example, \cite[Theorem~15.17]{ebud}.\footnote{
This we call a {\em Gr\"obner degeneration\/}.
The name ``Gr\"obner'' is attached because it seems somehow to evoke the right
picture.   In any case we are not the first to use this terminology.}
Further, $\tilde{V}$ is irreducible;  the
morphism~$\tilde{V}\to\mathbb{A}^1$ and all 
its base extensions  
are of relative dimension $\dim{Y}$ in the sense of \cite[Appendix B.2.5]{fulton:it}.
In particular, the scheme $M$ has pure dimension $\dim{Y}$ (meaning all its
irreducible components have dimension~$\dim{Y}$).   

Let $\mathbb{A}^1\subseteq\mathbb{P}^1$ be a standard open imbedding.
Let $U\times\mathbb{A}^1\subseteq U\times\mathbb{P}^1$ be the induced open
imbedding and ${\tilde{V}^\textrm{cl}}$ the closure of $\tilde{V}$ in $U\times\mathbb{P}^1$.    Then ${\tilde{V}^\textrm{cl}}$ is irreducible.
The morphism ${\tilde{V}^\textrm{cl}}\to
\mathbb{P}^1$ being dominant, it and all its base extensions  
are flat of relative dimension~$\dim{Y}$.

Applying the functor $E_n(\cdot)$ to the closed imbedding
$\tilde{V}^\textrm{cl}\hookrightarrow U\times \mathbb{P}^1$,   we get
a closed imbedding 
$E_n(\tilde{V}^\textrm{cl})\hookrightarrow E_n(U\times \mathbb{P}^1)$.
Since $T$ acts trivially (by definition) on~$\mathbb{P}^1$,  we have
$E_n(U\times\mathbb{P}^1)=E_n(U)\times\mathbb{P}^1$.   
Composing the inclusion
$E_n(\tilde{V}^\textrm{cl})\hookrightarrow E_n(U)\times \mathbb{P}^1$
with the second projection to~$\mathbb{P}^1$,  we get
a dominant (and hence flat) map $E_n(\tilde{V}^\textrm{cl})\to\mathbb{P}^1$;
composing it with the first projection,  we get
a proper map $E_n(\tilde{V}^\textrm{cl})\to E_n(U)$.

Since the flat pull-back of a class rationally equivalent to zero remains
rationally equivalent to zero~\cite[Theorem~1.7]{fulton:yt}, 
the fibers over two points in $\mathbb{P}^1$
are rationally equivalent in $E_n(\tilde{V}^\textrm{cl})$.  Since
rational equivalence is preserved by proper 
push-forward~\cite[Theorem~1.4]{fulton:yt},   the images
in $E_n(U)$ of these fibers are rationally equivalent.  
From the way $\tilde{V}^\textrm{cl}$ is constructed,   $E_n(Y\cap U)$ and $E_n(M)$ are respectively
the generic and special fibers in $\tilde{V}^\textrm{cl}$ over $\mathbb{P}^1$.
Thus the classes of $E_n(Y\cap U)$ and $E_n(M)$ are rationally equivalent
in $E_n(U)$.

It now remains only to show that $[E_n(M)]$ equals the right hand side 
of~(\ref{etgrob}).   Setting $M_i'$ to be the scheme defined
by~$\mathfrak{p}_i$, we have, by definition,
$[E_n(M)]=\sum_i \ell_i [E_n(M_i')]$.    The normal bundle to the smooth
subvariety $E_n(M_i')$ of $E_n(U)$ is the pull-back, via the projection, of
the bundle on $B_n$ associated to the $T$-module $U/M_i'\simeq
\langle X_1^{a_{1i}},\ldots,X_d^{a_{di}}\rangle^\textrm{dual}$ whose top Chern class
clearly is
$\chi_1^{a_{1i}}\cdots \chi_d^{a_{di}}$.    Hence $[E_n(M_i')]=
\chi_1^{a_{1i}}\cdots \chi_d^{a_{di}}$, completing the proof.\eproof}

\section{Towards an equivariant Pieri formula}\label{spieri}
Recall,  from \cite[Eq.~(9), p.146]{fulton:yt} for example,
that the classical Pieri formula gives a beautiful expression,
as an integral linear combination of general Schubert classes,
for the product of a special Schubert class with a general
Schubert class in the ordinary cohomology ring of the Grassmannian.
It seems like there ought to be a similarly beautiful
closed-form equivariant version that specializes
to the ordinary one.    Unfortunately,  this we do not yet have.

All we want to do in this section is record an observation
(see Proposition below) that is a formal consequence of the following basic
and well-known facts:  the injection of Equation~(\ref{einjection}); 
the restriction to a $T$-fixed point of an (equivariant) Schubert
class vanishes if the fixed point is not contained in the 
Schubert variety;   and the degree of a Schubert
class equals the codimension of the Schubert variety.
Being a formal consequence,   the observation holds for any
generalized flag variety,  not just the Grassmannian.  The point
to note is that the right hand side of (\ref{epieri}) is in terms
of restrictions,    which, thanks to Theorem~\ref{trestrdet}, we know
how to compute in the case of Grassmannians.

Let $G$ be a complex semisimple algebraic group and $Q$ a parabolic
subgroup.  Let~$T$ be a maximal torus and $B$ a Borel subgroup of $G$ such that
$T\subseteq B\subseteq Q$.   Let~$W$ denote the Weyl group of $G$ with respect
to $T$, and $W_Q$ the Weyl group of (the Levi part of) $Q$ with respect
to $T$.    The Schubert varieties in $G/Q$ are by definition the
$B$-orbit closures for the action of $B$ on $G/Q$ by left multiplication.
They are naturally indexed by $W/W_Q$.   We use $u$, $v$, $w$, \ldots\ to
denote elements of $W/W_Q$;  $X(u)$, $X(v)$, $X(w)$, \ldots\ denote the
corresponding Schubert varieties;  $[X(u)]$, $[X(v)]$, $[X(w)]$, \ldots\
denote the corresponding equivariant Schubert classes. 

The partial order on Schubert varieties
by inclusion induces a partial order, denoted $\leq\,$, on $W/W_Q$.
Let $c_{uv}^w$ be the structure constants of the multiplication 
of the equivariant integral cohomology ring of $G/Q$ with the respect to
the basis of Schubert classes:
\[ [X(u)]\cdot [X(v)] = \sum\limits_w c_{uv}^w [X(w)].	\]
The proof of the following proposition is a straightforward induction argument and so we omit it.
\bprop\label{ppieri}
\begin{enumerate} 
\item  $c_{uv}^w=0$ unless $w\leq u$, $w\leq v$, and the codimensions 
are such that
$\codim{X(u)}+\codim{X(v)}\geq \codim{X(w)}$.
\item If $w\leq u$ and $w\leq v$, then 
\begin{equation}\label{epieri}
c_{uv}^w =\sum\limits_{w=y_0<\ldots<y_k} (-1)^k \cdot
\frac{[X(u)]|_{y_k} [X(v)]|_{y_k}}{[X(y_k)]|_{y_k}}\cdot
\frac{[X(y_{k})]|_{y_{k-1}}}{[X(y_{k-1})]|_{y_{k-1}}}\cdot\,\cdots\,\cdot
\frac{[X(y_{1})]|_{y_{0}}}{[X(y_{0})]|_{y_{0}}}
\end{equation}
where the sum is over all chains $w=y_0<\ldots<y_k$ with
$y_k\leq u$ and $y_k\leq v$;   $[X(u)]|_{y_k}$ denotes the restriction
of the Schubert class $[X(u)]$ to the $T$-fixed point indexed by $y_k$.\hfill$\Box$
\end{enumerate}
\eprop

\bibliographystyle{amsplain}

\end{document}